\documentclass[a4paper]{article}
\usepackage[english]{babel}
\usepackage{pgfplots}\pgfplotsset{compat=1.18}
\usepackage{float}
\usepackage[T1]{fontenc}
\usepackage[utf8]{inputenc}
\usepackage{amsfonts}
\usepackage[affil-it]{authblk}
\usepackage{algpseudocode}
\usepackage{paralist}
\usepackage{caption}
\usepackage{appendix}
\usepackage{verbatim}
\usepackage{wrapfig}
\usepackage{tikz-cd} 
\usepackage{tikz}
\usetikzlibrary{plotmarks}
\usepackage{caption}
\usepackage{amsmath}
\usepackage[arrow, matrix, curve]{xy}
\usepackage{mathtools}

\usepackage[x11names]{xcolor}
\usepackage{amssymb}
\usepackage{wasysym}
\usepackage{listings}
\usepackage[gen]{eurosym}

\usepackage{amssymb,bbm}
\usepackage{stmaryrd}
\usepackage{graphicx}
\usepackage{mathrsfs}
\usepackage{pst-func}
\newcommand{\dd}{\mathrm{d}}

\DeclareMathOperator{\Sym}{\R^{\textit{n}\times \textit{n}}_{\mathrm{Sym}}}
\DeclareMathOperator{\Symtwo}{\R^{2\times 2}_{\mathrm{Sym}}}
\DeclareMathOperator{\Rnn}{\R^{\textit{n}\times \textit{n}}}

\DeclareMathOperator{\Snegativthree}{\mathcal{S}_-^{3\times 3}}
\DeclareMathOperator{\tr}{\mathrm{tr}}
\DeclareMathOperator{\dev}{\mathrm{dev}}
\DeclareMathOperator{\Symthree}{\R_{\mathrm{Sym}}^{3\times 3}}

\DeclareMathOperator{\id}{\mathrm{id}}
\DeclareMathOperator{\R}{\mathbb{R}}
\DeclareMathOperator{\N}{\mathbb{N}}
\DeclareMathOperator{\relu}{\mathrm{ReLU}}
\DeclareMathOperator{\diag}{\mathrm{diag}}
\DeclareMathOperator{\cof}{\mathrm{cof}}

\usepackage{enumitem}
\usepackage{tabularx}
\usepackage{amsthm}

\usepackage{tikz}
\usetikzlibrary{decorations, decorations.text}

\newtheorem{Theorem}{Theorem}[section]

\newtheorem{lemma}[Theorem]{Lemma}
\newtheorem{Bsp}[Theorem]{Example}
\newtheorem{Bem}[Theorem]{Remark}
\newtheorem{Proposition}[Theorem]{Proposition}
\newtheorem{Definition}[Theorem]{Definition}

\usepackage{fancyhdr}
\usepackage{datetime}
\renewcommand{\sectionmark}[1]{\markboth{~#1}{}}

\fancypagestyle{main}{
\fancyhf{}
\fancyhead[LE]{\footnotesize\nouppercase{\sectionmark}}
\fancyhead[RO]{\footnotesize \nouppercase{\sectionmark}}
\fancyfoot[LE,RO]{\thepage}

}
\renewcommand\subsectionmark[1]{}

\usepackage{csquotes}
\usepackage[backend = biber, sorting=nyt, style =numeric-comp, giveninits=true,
  maxbibnames=99]{biblatex}

\usepackage[breaklinks]{hyperref}

\begin{document}

\title{\bfseries Universal neural network approximation of symmetric polyconvex energies and envelopes}

\author[1]{{Florian Gammer}}
\date{23\textsuperscript{rd} October 2026}

\affil[1]{\footnotesize Faculty of Mathematics, University of Regensburg, Universitätsstraße 31
93053 Regensburg, Germany}

\maketitle

\begin{abstract} 
Symmetric polyconvexity provides a sufficient condition for weak lower semicontinuity in small-strain elasticity and a tractable lower bound for the relaxation of non-(quasi)convex energies. 
Since classical algorithms for computing such functions can be computationally expensive, neural network approximations have recently emerged as an alternative.
A central difficulty with establishing such neural network based approaches for symmetric polyconvexity is the additional monotonicity property. As we show, the corresponding weight restrictions on ReLU-based input-convex networks do not yield universal approximation of the lifted representers. 
Here, we introduce symmetric polyconvex neural networks (SPCNNs) for constitutive modeling and approximation of symmetric polyconvex envelopes. The deep architecture combines convexity with the monotonicity required of a representer in the space of minors. We prove a universal approximation theorem for symmetric polyconvex energies, without requiring isotropy, and quantify the number of affine pieces attainable by the lifted network and derive reduced architectures for isotropic energies. Numerical examples illustrate stress-based regression and approximation of the symmetric polyconvex envelope of a double-well potential. Symmetric polyconvexity is guaranteed by construction, while the accuracy of the fitted energies and envelopes is assessed numerically.
\end{abstract}

\smallskip


\section{Introduction}

In solid mechanics, physical phenomena are mathematically described by an interplay of energy and dissipation. 
In variational models of microstructure formation, lack of quasiconvexity can reflect minimizing sequences with increasingly fine oscillations, leading to the concept of relaxation of the energy functional, which corresponds to computing the quasiconvex envelope \cite{Daco08,benesova2017weak,Cia88,morrey1952quasi}. In the context of small-strain elastoplasticity, the strains and stresses are symmetric matrices, and relaxation corresponds to computing the \emph{symmetric} quasiconvex envelope \cite{barroso2000relaxation,fonseca1999cala,ebobisse2000lower,rindler2011lower,zhang2002equality,peigney2013energy}. Although calculating this envelope directly is generally unfeasible, the common approach to approximate it, is to compute an upper bound (the symmetric rank-one convex envelope) as well as a lower bound (the symmetric polyconvex envelope), and then show that these match (e.g., \cite{chenchiah2008relaxation,Kohn91,KoStrII,grabovsky2019rank,grabovsky1996bounds,grabovsky2019explicit,grabovsky1995microstructures,martin2017rank,GOVINDJEE2003I,fechte2016lamination,antimonov2016phase,KristensenGuerra2022}). 

Originating in Ball's finite elasticity theory \cite{Ball1976}, polyconvexity allows one to construct quasi-convex energy densities respecting physical constraints. Unlike standard polyconvexity, symmetric polyconvexity requires functions to be convex in their first- and second-order minors and non-increasing with respect to second-order minors \cite{Kreisbeck19,Merabet24}.

Although in theory symmetric polyconvexification is more accessible conceptually than symmetric quasiconvexification, computing it remains a demanding optimization problem (cf.~\cite{bartels2005reliable}).
A neural network approach offers a reusable approximation that can be evaluated at any strain after training.
For constitutive model identification, the objective is to find an energy within a prescribed class of symmetric polyconvex functions. For relaxation, the objective is to approximate the largest symmetric polyconvex function lying below a given non (poly-)convex energy. This envelope provides a lower bound for the symmetric quasiconvex envelope and can help identify the relaxed energy when combined with matching upper bounds.

Here, we develop a neural network architecture for these two tasks. Building on input-convex neural networks \cite{amos17b} and the characterization of symmetric polyconvexity \cite{Kreisbeck19}, we enforce convexity and monotonicity of a lifted representer.

Our main theoretical result is universal approximation of symmetric polyconvex energies on compact subsets of the space of symmetric matrices, without an isotropy assumption. Obtaining this result requires more than imposing cofactor monotonicity through weight restrictions on a $\relu$-based ICNN. We show that the resulting restricted networks satisfy an additional order inequality and construct a symmetric polyconvex energy for which every finite convex representer violates this inequality. Thus, these networks cannot approximate all lifted representers locally uniformly. Our architecture overcomes this obstruction by replacing $\relu$-activations with pairwise maxima, while retaining the required convexity and monotonicity. We further quantify the expressivity per parameter of the lifted output, and derive reduced architectures for isotropic and two-dimensional energies.

As a quick outline of this paper, Section 2 recalls symmetric polyconvexity and its relation to relaxation, reviews related neural network approaches, and discusses the expressivity of ICNNs. Section 3 establishes the insufficiency of restricted ReLU networks, introduces the SPCNN architecture, and proves its universality and expressivity properties. Section 4 describes stress-based constitutive regression and the training objective for envelope approximation. Section 5 presents numerical results for both applications.

\section{Preliminaries}

We are essentially concerned with the problem of minimizing an energy $\mathbb{E}[u]:=\int_\Omega F(\varepsilon(u))\dd x$,
 where $\varepsilon(u)$ denotes the \emph{symmetrized} gradient, $\Omega\subset \R^n$ is a Lipschitz domain, $n\in \{2,3\}$, $F\colon \Sym \to \R$ is an energy density, 
$u$ is an element of the Sobolev space $ W_g^{1,p}(\Omega,\R^n):=g+W_0^{1,p}(\Omega,\R^n)$, $p> 1$, where $g\in W^{1,p}(\Omega,\R^n)$. Let us write 
$$\pi_S\colon \R^{n\times n}\to \Sym, \ A\mapsto \frac{1}{2}(A+A^T)$$
for the orthogonal projection. Then we can express
$$\int_\Omega F(\varepsilon(u))\dd x =\int_\Omega F\circ \pi_S(\nabla u)\dd x.  $$
Using the direct method in the calculus of variations, sufficient requirements for the existence of a minimizer are coercivity and (depending on the growth of $F$) quasiconvexity of $F\circ \pi_S$, which is the same as $F$ being symmetric quasiconvex.

\subsection{Symmetric polyconvexity, isotropy, and envelopes}
 
Analogously to  symmetric quasiconvexity, we can define symmetric polyconvexity using the projection. Throughout, $\Sym$ denotes the space of symmetric matrices, lowercase Greek letters $\varepsilon,\xi,\eta$ arbitrary elements of this space, and $\varepsilon:\xi$ the Frobenius scalar product. 
\begin{Definition}
    A function $F\colon \Sym\to \R$ is called symmetric polyconvex if $F\circ \pi_S$ is polyconvex, i.e., there exists a convex function 
    $$\hat{F}\colon \R^{n\times n} \times \R^{\binom{n}{2}\times \binom{n}{2}} \times ... \times \R^{\binom{n}{n-1}\times \binom{n}{n-1}}\times \R^{\binom{n}{n}\times \binom{n}{n}}\to \R  $$
    such that 
    $F\circ \pi_S(A)=\hat{F}(T(A))$ for all $A\in \Rnn$, where $T(A)$ is the vector of minors of the matrix $A$. 
\end{Definition}

When we are considering problems of the symmetric form, it is convenient to work entirely in the space of symmetric matrices, instead of taking the whole space into account at first and then reducing to the symmetric case later on. A self-contained characterization of symmetric polyconvexity 
 was proved by Boussaid, Kreisbeck and Schlömerkemper \cite{Kreisbeck19}.
 We denote the set of symmetric negative semidefinite matrices by
 $\Snegativthree$, and $\partial_2$ is the subdifferential with respect to the second variable. Then it holds that 
\begin{Theorem} \label{Kreisbeck Theorem} 
    \begin{enumerate}
        \item[(i)] a function $F\colon \Symtwo\to \R$ is symmetric polyconvex if and only if there is a convex function $G\colon \Symtwo \times \R\to \R$ that is nonincreasing in the last argument, such that $F(\varepsilon)=G(\varepsilon,\det \varepsilon)$.
        \item[(ii)] a function $F\colon \Symthree\to \R$ is symmetric polyconvex if and only if there is a convex function $G\colon \Symthree \times \Symthree\to \R$ that is non-increasing in the second argument, i.e., admits $\partial_2 G(\varepsilon,\eta)\subset \Snegativthree$ for all $\varepsilon,\eta\in \Symthree$,
        such that $F(\varepsilon)=G(\varepsilon,\cof  \varepsilon)$.
     \end{enumerate}
\end{Theorem}
While symmetric polyconvexity of $F$ (accompanied with some growth conditions) gives rise to the existence of minimizers, in solid mechanical practice, many arising potentials are neither symmetric polyconvex nor symmetric quasiconvex. Such models are derived from microscopic considerations and, crucially, the lack of symmetric quasiconvexity, and the associated potential non-existence of minimizers, is a necessary feature for modeling, e.g., the formation of microstructure. 
In more abstract terms, instead of a minimizer, only an oscillating minimizing sequence may exist, yielding a macroscopic behavior that corresponds to a microscopic mixture of phases. This mixture behavior can be quantified by the symmetric quasi- or polyconvex envelopes.

We continue our discussion by investigating the definition of the symmetric polyconvex envelope of a function $F\colon \Sym\to \R$.
We remark that our functions are finite-valued and assumed to be bounded from below by polyaffine functions, thus the (symmetric poly-) convex envelope is the same as the lower semicontinuous (symmetric poly-) convex envelope (\cite[Thm.~2.43]{Daco08}), and we will not explicitly mention the lower semicontinuity property of the envelopes in what follows.

\begin{Definition}\label{Def: sym polyconvex envelope}
    Let $F\colon \Sym\to \R$ be bounded from below. Then the symmetric polyconvex envelope is defined as
    \begin{align*}
F^{\mathrm{spc}} \colon \Sym &\to \R\\
    \varepsilon&\mapsto \sup\{F'(\varepsilon)\mid F'\colon \Sym\to \R \ \text{is symmetric polyconvex}\\
    &\ \ \ \ \ \ \ \ \ \ \ \ \ \ \ \ \ \ \ \ \ \ \ \ \ \ \text{and}\ F'(\xi)\le F(\xi) \ \text{for all}\ \xi\in \Sym \}.
\end{align*}
\end{Definition}
Classical numerical methods typically use more explicit representations, such as Carathéodory formulas (cf.~\cite[Thm.~6.8]{Daco08},\cite{bartels2005reliable,neumeier2024computational,GamKre2026}).
This abstract form is, however, the foundation of our neural network relaxation procedure. To reduce computations and make simpler plots, we will constrain our later benchmark examples to the isotropic case. A function is isotropic if for all $\xi\in \Sym$ and for all $Q\in \mathrm{O}(n)$ it holds that 
$$F(Q^T\xi Q)=F(\xi).$$
In this case, the complexity reduces, since it suffices to calculate the above envelope on diagonal matrices.

\begin{Theorem}[\cite{GamKre2026}]\label{Thm: Isotropic Representation}
 For $\varepsilon\in \Symthree$ let $\Lambda(\varepsilon)$ be the vector of eigenvalues of $\varepsilon$ in increasing order. 
 Let $F\colon \Symthree\to \R$ be isotropic. Let $f\colon \R^3\to \R$ be defined by
 $$f(x_1,x_2,x_3):=F(\diag(x_1,x_2,x_3)).$$
 Then $$F=f\circ \Lambda$$
 and $f$ is permutation invariant. 
 Moreover, $F$ is symmetric polyconvex if and only if there is a convex function $g\colon \R^3\times \R^3\to \R$, which is jointly permutation invariant,
    and non-increasing w.r.t.~the second variable, such that for all $x\in \R^3$:
    $f(x)=g(x,\cof x)$, where $\cof(x_1,x_2,x_3)=(x_2x_3,x_1x_3,x_1x_2)$. 
 \end{Theorem}

\subsection{Overview of the literature on polyconvex neural networks}

Traditional numerical algorithms for finding solutions to PDEs and differential inclusions, such as the finite element method, are widely used in practice due to their robust error bounds, interpretability, and relative simplicity. However, they come with a set of limitations, such as difficulty in handling noisy data, the curse of dimensionality, the complexity of mesh generation, and the mesh-based representation of their output. Neural networks offer an alternative parameterization whose approximation efficiency depends on the target function, architecture, and loss function.

Fundamentally, these networks solve a regression problem within a large class of admissible functions. For an introduction and overview of such general feedforward networks, see \cite{jentzen2025mathematicalintroductiondeeplearning}. Depending on the activation function, the class of functions is chosen to be dense in the function space one aims to approximate. In machine learning, this is referred to as the universal approximation property. The foundational proofs of this property were given in \cite{Cybenko1989} (using the Hahn--Banach and Riesz--Markov Theorem), \cite{HORNIK1989359} (using the Stone-Weierstrass Theorem), and further specified for instance in \cite{wang2004general,Hanin2019, hornik1991approximation}.
While neural networks have been highly successful in natural language and image processing, their application in physical modeling has, until recently, been limited. This is because physical systems adhere to precise mathematical laws that a standard network's output might violate, especially outside the training domain. Furthermore, there is often insufficient physically relevant data to train vast neural networks. In recent years, these challenges have been simultaneously addressed through the development of Physics-Informed Neural Networks (PINNs) (cf.~\cite{stavroulakis2025physics}) and structurally constrained architectures (\cite{karniadakis2021physics}). While standard PINNs enforce physical laws via soft constraints in the loss function, heavily constrained architectures ensure the output inherently satisfies desired properties, such as Lipschitz continuity \cite{anil2019sorting,tanielian2021approximating}, convexity \cite{amos17b}, monotonicity \cite{daniels2010monotone}, or certain symmetries. 

Abstractly speaking, these constraints restrict the class of approximable functions to those satisfying the prescribed physical and mathematical properties.

Data-driven approaches in computational mechanics range from methods that use experimental material data directly to methods that learn constitutive relations (\cite{kirchdoerfer2016data, BARTEL2023116138, as2022mechanics}). In energy-based approaches, a stored energy density is learned, and stresses are obtained by differentiation.

Polyconvexity provides a useful sufficient condition for quasiconvexity and, together with suitable coercivity and growth assumptions, ensures existence of minimizers. This motivates neural network architectures that enforce polyconvexity by construction. Existing approaches use deformation-based inputs, including invariants of the right Cauchy–Green tensor and, for anisotropic materials, additional structural invariants (\cite{chen2022polyconvex,fuhg2022learning, kalina2022automated, LINDEN2023105363,klein2022polyconvex}), as well as principal stretches and their combinations (\cite{tepole2025polyconvex}).
Most recently, \cite{geuken2025input} utilized the characterization of isotropic polyconvex functions \cite{wiedemann2026characterization} to derive a streamlined neural network architecture that inherently produces such functions by reducing the problem to principal stretches and their products.

Our approach in this paper is related, but distinct in two crucial ways. First, our network concerns symmetric polyconvexity, and therefore small-strain elasticity. 
Working in this framework, dependence on the symmetric gradient ensures invariance under infinitesimal rigid motions, while finite-strain conditions such as energy blow-up under compression to zero volume lie outside the small-strain setting in the first place.
Second, our neural network goes beyond simply lifting an ICNN to the space of minors, since we also need to account for monotonicity in the second order minors, while preserving the universal approximation property.

\subsection{On the expressivity and universality of ICNNs}

Common among most of the aforementioned approaches to polyconvex neural networks is their foundation, the so-called input convex neural networks (ICNN) introduced in \cite{amos17b}. These networks are a special type of feedforward network where the output function is automatically convex in the inputs. 
In this subsection, we give a quick introduction to these networks and compare their expressivity to previous architectures. 

The simplest neural network architecture that guarantees convexity in its outputs and has the universal approximation property is the single-layer max-out network as a special case of \cite{goodfellow2013maxout}. This is nothing more than having $n\in \N$ affine functions $L_i$ and setting the network output to $\max_{i=1,...,n}L_i(y)=\mathcal{F}(y,\Theta)$, where $\Theta$ is the collection of the parameters, i.e., the slopes and intercepts of the $L_i$. Essentially, this is the same as finding a piecewise affine and convex approximation to a convex function. The general problem of fitting arbitrary convex functions by piecewise affine ones numerically has been addressed in numerous papers, such as \cite{cox1971algorithm, magnani2009convex, gavrilovic1975optimal}. 

The downside of these simple architectures is their lack of expressivity due to their shallowness. Of course, every convex function on a compact convex set can be approximated by these networks, but this may require a large number of parameters. 

The ICNNs introduced in \cite{amos17b} solve this issue by presenting a deep neural network architecture as follows. 
Let $d$ be the space dimension of the inputs, let $y\in \R^d$ be an input vector, let $k\in \N$, $n_1,...,n_k\in \N$, $n_0=d,n_{k+1}=1$ and let $W_i^{(z)}\in \R^{n_{i+1}\times n_{i}} $, $W_i^{(y)}\in \R^{n_{i+1}\times d} $, $b_i\in \R^{n_{i+1}}$ be parameters for $i=0,...,k$, and $W_0^{(z)}=0$. Furthermore, let 
$g_i\colon \R\to \R$ be activation functions for $i=0,...,k-1$ that are applied componentwise, and let $g_k=\id_{\R}$. 
Then we iteratively set
\begin{align*}
    z_0=0, \quad z_{i+1}& = g_i\bigl( W_i^{(z)}z_i + W_i^{(y)}y+b_i\bigr) \text{ and } \mathcal{F}(y,\Theta) = z_{k+1}.
\end{align*}
Here, $\Theta$ is again the collection of network parameters and $\mathcal{F}(\cdot,\Theta)$ the output function of the network with parameters $\Theta$. The first important difference to general fully connected feedforward neural networks is the feature of direct pass-through layers. The second difference is the restriction of the parameters and activation functions. 
\begin{Proposition}[\cite{amos17b} ]\label{Prop: ICNN}
Let $W_i^{(z)}\in [0,\infty)^{n_{i+1}\times n_{i}}$ and let $g_i$ be convex and non-decreasing for all $i=0,...,k$ and $g_k=\mathrm{id}_{\R}$. Then the output function $\mathcal{F}(\cdot,\Theta)$ of the network described above is convex in the inputs. Under these conditions, the network is called an \emph{ICNN}.
\end{Proposition}
When we choose the $g_i$ to be the $\mathrm{ReLU}$ function (i.e., $\max\{x,0\}$), then an ICNN of width $1$ and depth $n$ can represent a single layer max-out network with $n$ neurons directly (cf., \cite{chen2019optimalcontrolneuralnetworks} Thm.~2 (1)). This is also one of the reasons why ICNNs admit the universal approximation property.
\begin{Theorem}[\cite{huang2021convex,chen2019optimalcontrolneuralnetworks}]
    Let $f\colon [0,1]^d\to \R$ be continuous, convex and let $\varepsilon>0$. Then there is an \emph{ICNN} $\mathcal{F}(\cdot,\Theta)$ with $\mathrm{ReLU}$-activation, such that
    $$\sup_{y\in [0,1]^d}\vert f(y)-\mathcal{F}(y,\Theta)\vert <\varepsilon.$$
\end{Theorem}
Moreover, an ICNN with only one hidden layer ($k=1$) and width $N$ can represent a single layer max-out network with $\mathcal{O}(N^d)$ neurons. 
To our knowledge, a precise characterization accounting for these geometric activation constraints has not been explicitly established in the literature\footnote{In \cite[Thm.~2(2)]{chen2019optimalcontrolneuralnetworks}, a related bound of $2^N$ affine pieces is stated, which overestimates the realizable capacity, since once an affine component is activated past a hyperplane $L_i=0$, it remains active on that side.}.

\begin{Theorem}[Expressivity of single layer ICNNs]\label{Thm: Expressivity of ICNN}
 There is a single layer \emph{ICNN} $\mathcal{F}(\cdot,\Theta)$ with  $\mathrm{ReLU}$-activation of width $N\ge d$ and $C(N,d)$ distinct affine functions $L_i$, $i=1,...,C(N,d)$, such that $\mathcal{F}(\cdot,\Theta)=\max_{i=1,...,C(N,d)}L_i$. Here,   $C(N,d)$ is the maximal number of regions in $\R^d$ that can be separated by $N$ hyperplanes:
 $$C(N,d):=\sum_{j=0}^d\binom{N}{j}=\frac{N^d}{d!}+\mathcal{O} (N^{d-1}).$$
\end{Theorem}
\begin{proof}
    The single layer ICNN can be written as 
    $$\mathcal{F}(y,\Theta)=\sum_{i=1}^N\relu (L_i(y)),$$
    for affine functions $L_i$, where we set the pass-through layers to zero.
We choose $L_1, \ldots, L_N$ so that all subset sums $\sum_{i \in I} L_i$ are pairwise distinct. This is possible because forbidden coincidences lie within a finite collection of proper linear subspaces, leaving a dense open set of admissible choices.
    
Each hyperplane $\{L_i=0\}$ separates $\R^d$ into two distinct regions. Crucially, $N$ hyperplanes separate $\R^d$ into maximally $C(N,d)$ regions (we can define $C(N,d)$ like this), and we may choose the $L_i$ such that the maximal number of regions is attained. This means that $\mathcal{F}(\cdot,\Theta)$ is continuous, convex, piecewise affine on every separated region, and pairwise distinct on each region. In particular, the corresponding single layer max-out network requires $C(N,d)$ many affine functions. For a derivation of the closed form for $C(N,d)$, see \cite{Buck1943PartitionOS}. 
\end{proof}

\section{Symmetric polyconvex neural networks}

In this section, we propose an efficient and universal architecture for an ICNN whose output is symmetric polyconvex in its inputs. The first observation hereby is to use the characterization of symmetric polyconvex functions in Thm.~\ref{Kreisbeck Theorem}. Essentially, these are higher dimensional convex functions of the first and second order minors of the input matrix. For regular polyconvexity, one can build an ICNN approximation for these lifted functions, cf., \cite{geuken2025input, klein2022polyconvex,tepole2025polyconvex}.  
In contrast, symmetric polyconvexity additionally requires the lifted output function to decrease in the second order minors. General monotonic neural networks have been considered in the literature for instance in \cite{wehenkel2019unconstrained,daniels2010monotone}, and monotonic ICNNs for polyconvex envelopes have been considered recently in \cite{neumeier2026compression}. However, this parameter restriction is not sufficient for universal approximation of the lifted convex representers.

\subsection{An obstruction to lifted approximation by restricted ReLU networks}

A symmetric polyconvex network could potentially be realized by setting the parameters corresponding to the minors in the lifted ICNN to be negative semidefinite. We show that this approach cannot approximate every convex, cofactor-nonincreasing lifted representer, even when arbitrary depth and width are allowed. 
The reason is that such $\relu$-networks inherently satisfy a particular order relation. 
\begin{lemma}[ReLU order]\label{lem: ICNN with non-positive pass-through layer weights}
    Let $\mathcal{F}\colon \R^2\to \R$ be the output function of an $\mathrm{ICNN}$ of arbitrary depth and width, as defined in Prop.~\ref{Prop: ICNN}, with non-positive pass-through layer weights and $\relu$-activation, and non-positive coefficients
for both input variables in every input and pass-through
connection, including those to the output layer.
    Then $\mathcal{F}$ is convex, non-increasing, and satisfies the inequality
    \begin{align}\label{eq: Relu inequality} \tag{$\ast$}
        \mathcal{F}(x_0,y_0)+\mathcal{F}(x_1,y_1)\ge \mathcal{F}(x_1,y_0)+\mathcal{F}(x_0,y_1)
    \end{align}
    for all real $x_0\le x_1$, $y_0\le y_1$.
\end{lemma}
\begin{proof}
    The proof follows by induction. If every component of $z_i$ has these properties,
then so does every preactivation in the next layer, since
it is a nonnegative linear combination of these components
plus an affine function with non-positive coefficients
in the inputs. The induction starts with $z_0=0$.
It remains to show that applying $\relu$ preserves these
properties. For convexity and monotonicity, this is clear. Let $g\colon\R^2\to\R$ denote a preactivation
having all three properties.
  Let us write
    $$A=g(x_0,y_0),\ B=g(x_1,y_0), \ C=g(x_0,y_1), \ D=g(x_1,y_1).$$
    By assumption, we have
    $$A\ge B, \ C\ge D, \ A+D\ge B+C.$$
If $B,C\le 0$, the claim follows from the nonnegativity of $\relu$. If $B\le 0$ and $C>0$, then 
\begin{align*}
    \relu(A)+\relu(D) \ge A\ge C=\relu(C)+\relu(B),
\end{align*}and analogously for switched roles of $B,C$. Finally, suppose $B,C>0$, then 
\begin{align*}
    \relu(A)+\relu(D)\ge A+D\ge B+C= \relu(C)+\relu(B).
\end{align*}
\end{proof}

For the following symmetric polyconvex function, every finite convex representer violates the $\relu$-order inequality and therefore lies outside the locally uniform closure of the restricted lifted networks.

\begin{Bsp}
  Consider the two affine (in the minors) functions 
$$L_1(\varepsilon,\cof\varepsilon)=\diag(0,-2,-1):\cof \varepsilon, \quad L_2(\varepsilon,\cof\varepsilon)=\diag(0,-1,-2):\cof \varepsilon.$$
For the symmetric polyconvex function $F(\varepsilon):=\max(L_1(\varepsilon,\cof\varepsilon),L_2(\varepsilon,\cof\varepsilon))$, no finite convex representer can be approximated locally uniformly by parameter-restricted $\relu$-$\mathrm{ICNN}$s.
We note that any convex representer $G\colon \Symthree\times \Symthree\to \R$ of $F$ needs to be constant w.r.t.~the first variable. To see this, take a rank-one direction $a\otimes a$, $a\in \R^3$. Then $\cof(a\otimes a)=0$, hence $F(ta\otimes a)=G(ta\otimes a,0)=0$ for all $t\in \R$. For an arbitrary $(\xi,\eta)$ choose $(\alpha,\beta)\in\partial G(\xi,\eta).$ 
The subgradient inequality gives
$$0\ge
G(\xi,\eta)+\alpha:(ta\otimes a-\xi)-\beta:\eta
\quad \text{for every } t\in \R.$$
We infer $\alpha:(a\otimes a)=0$ for every $a\in \R^3$, hence $\alpha=0$.

Furthermore, the representer is uniquely determined on the positive semidefinite matrices by continuity, since $\eta=\cof (\sqrt{\det \eta}\,\eta^{-1})$ for positive definite $\eta$.
Thus, every representer $G$ satisfies for all $u,v\ge0$
$$G\bigl(0,\diag(0,u,v)\bigr)
=\max\{-2u-v,-u-2v\}=:m(u,v).$$
For the latter function, assume an $\mathrm{ICNN}$ $\mathcal{F}$ as in lemma \ref{lem: ICNN with non-positive pass-through layer weights} exists such that
$$\delta=\sup_{u,v\in [0,1]}\vert \mathcal{F}(u,v)-m(u,v)\vert.$$
Then, using \eqref{eq: Relu inequality}, we see that
\begin{align*}
    0&\le  \mathcal{F}(0,0)+\mathcal{F}(1,1)-\mathcal{F}(1,0)-\mathcal{F}(0,1)\\
    &\le m(0,0)+m(1,1)-m(1,0)-m(0,1) +4\delta \\
    &= -1+4\delta,
\end{align*}hence $\delta\ge \frac{1}{4}$, and cannot be made arbitrarily small. 
\end{Bsp}

\begin{Bem}
    We note that the $\relu$-architecture still yields universal approximability in 2d due to the total ordering. Precisely, for any finite collection of affine functions of the form $L_i(\varepsilon,\det\varepsilon)=\alpha_i:\varepsilon+b_i\det\varepsilon+c_i$, where $b_i\le 0$, 
we relabel the affine functions so that $b_1\le\cdots\le b_n$. The result can be written iteratively as a chain of affine and $\mathrm{ReLU}$-functions analogous to Thm.~2.1 in \cite{chen2019optimalcontrolneuralnetworks}. Due to the ordering, the terms $L_i-L_{i+1}$ are non-increasing w.r.t.~the second variable. However, this rearrangement can cause problems in the training of the network, therefore the modified architecture of the next section is still preferable. 
\end{Bem}

\subsection{SPCNN architecture}

The architectural changes that we propose still use parameter restriction first, but in order to obtain universality, we require activation functions with certain properties. As we showed, the $\mathrm{ReLU}$-function is not sufficient, since it compares a neuron output only to the zero function. Instead, we compare the output of several neurons to each other. Similar approaches to the activation function have been carried out in the literature for deep max-out neural networks in \cite{goodfellow2013maxout} using general maxima, and more recently for Lipschitz continuous neural networks \cite{tanielian2021approximating, anil2019sorting, chernodub2017normpreservingorthogonalpermutationlinear} using groupsort activations. 

\begin{Definition}[SPCNN architecture in 3d]
    Let $\xi,\eta\in \Symthree$ be input matrices, let $k\in \N$, $n_1,...,n_k\in \N$, $n_0=n_{k+1}=1$, and let 
    \begin{align*}
        W_{i}^{(z)}&\in [0,\infty)^{2n_{i+1}\times n_{i}}, \quad b_i\in \R^{2n_{i+1}},\\
W_{i}^{(\xi)}&=\bigl(W_{i,j}^{(\xi)}\bigr)_{j=1}^{2n_{i+1}}\in \bigl(\Symthree\bigr)^{2n_{i+1}} ,  \quad W_{i}^{(\eta)}=\bigl(W_{i,j}^{(\eta)}\bigr)_{j=1}^{2n_{i+1}}\in \bigl(\Snegativthree\bigr)^{2n_{i+1}} 
    \end{align*}
be parameters for $i=0,...,k$, and $W_0^{(z)}=0$. We can group the $2n_{i+1}$ neurons into $n_{i+1}$ pairs of adjacent neurons (i.e., the first two, the second two, etc.). Furthermore, let 
$g_i\colon \R^2\to \R$ be activation functions that are convex and non-decreasing applied to two adjacent components for $i=0,...,k$.
We iteratively set
\begin{align*}
    z_0=0, \ z_{i+1}& = g_i\biggl( W_i^{(z)}z_i +\bigl( W_{i,j}^{(\xi)}:\xi+W_{i,j}^{(\eta)}:\eta\bigr)_{j=1}^{2n_{i+1}}+b_i\biggr),\ \mathcal{F}(\xi,\eta,\Theta) = z_{k+1}.
\end{align*}Then the set $\Theta$ of all parameters together with the iteration scheme is called a 3d symmetric polyconvex neural network (\emph{SPCNN}) with activations $(g_i)_{i=0}^k$, $\mathcal{F}_{\mathrm{lift}}(\cdot,\Theta)\colon \Symthree\times \Symthree\to \R, (\xi,\eta)\mapsto \mathcal{F}(\xi,\eta,\Theta)$ is called the lifted output function of the \emph{SPCNN} and $\mathcal{F}_{\mathrm{SPCNN}}(\cdot,\Theta)\colon \Symthree\to \R, \varepsilon\mapsto \mathcal{F}_{\mathrm{lift}}(\varepsilon,\cof\varepsilon,\Theta)$ is called the symmetric polyconvex output of the \emph{SPCNN}.
\end{Definition}

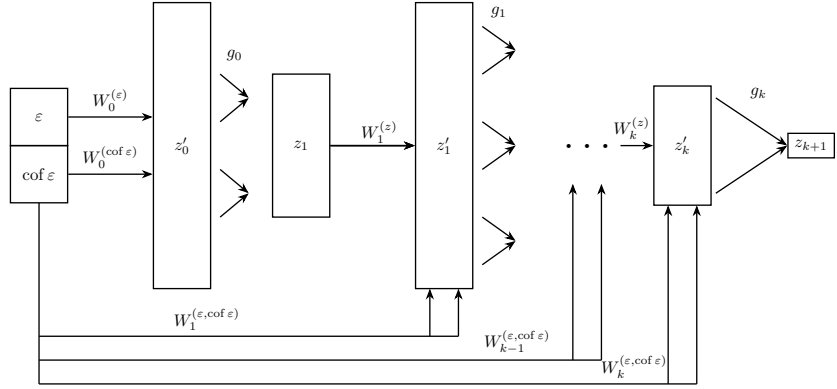
\begin{figure}[h]
\centering
\resizebox{0.9\textwidth}{!}{%
\begin{tikzpicture}[
    >=Stealth,
    layer/.style={draw, thick, minimum width=1.2cm, align=center},
    input/.style={draw, thick, minimum width=1.2cm, minimum height=1.2cm, align=center}
]
\node[input] (X) at (0, 0.6) {$\varepsilon$};
\node[input] (Y) at (0, -0.6) {$\cof\varepsilon$};

\node[layer, minimum height=6.0cm] (z0)  at (3.0, 0) {  $z_0'$};
\node[layer, minimum height=3.0cm] (z01) at (5.5, 0) {  $z_{1}$};
\node[layer, minimum height=6.0cm] (z1)  at (8.5, 0) {  $z_1'$};
\node (dots) at (11.5, 0) {\Huge $\ldots$};
\node[layer, minimum height=2.5cm] (zk1) at (13.5, 0) {  $z_{k}'$};
\node[draw, thick, inner sep=4pt] (zk) at (16.2, 0) {  $z_{k+1}$};

\draw[->, thick] (X.east) -- node[above] {  $W_0^{(\varepsilon)}$} (X.east -| z0.west);
\draw[->, thick] (Y.east) -- node[above] {  $W_0^{(\cof\varepsilon)}$} (Y.east -| z0.west);

\draw[->, thick] (3.8, 1.5) -- (4.4, 1.0);
\draw[->, thick] (3.8, 0.5) -- (4.4, 1.0);
\node at (4.1, 1.9) {  $g_0$};

\draw[->, thick] (3.8, -0.5) -- (4.4, -1.0);
\draw[->, thick] (3.8, -1.5) -- (4.4, -1.0);

\draw[->, thick] (z01.east |- 1.0,0) -- (z1.west |- 2.0,0) node[midway, above, xshift=4pt, yshift=-2pt] {  $W_1^{(z)}$};
\draw[->, thick] (z01.east |- 0,0) -- (z1.west |- 0,0);
\draw[->, thick] (z01.east |- -1.0,0) -- (z1.west |- -2.0,0);

\draw[->, thick] (9.3, 2.5) -- (10.0, 2.0);
\draw[->, thick] (9.3, 1.5) -- (10.0, 2.0);
\node at (9.65, 2.8) {  $g_1$};

\draw[->, thick] (9.3, 0.5) -- (10.0, 0.0);
\draw[->, thick] (9.3, -0.5) -- (10.0, 0.0);

\draw[->, thick] (9.3, -1.5) -- (10.0, -2.0);
\draw[->, thick] (9.3, -2.5) -- (10.0, -2.0);

\draw[->, thick] (12.2, 0) -- (zk1.west |- 0,0)node[midway, above,  xshift=-3pt] {  $W_k^{(z)}$};

\draw[->, thick] (14.2, 1.0) -- (zk.west |- 1.0,0);
\draw[->, thick] (14.2, -1.0) -- (zk.west |- 1.0,0);
\node at (15.1, 1.1) {  $g_k$};

\coordinate (bus_start) at (0.0, -1.2); 
\draw[thick] (bus_start) -- (0.0, -5);

\draw[thick] (0.0, -4) -- (8.8, -4) node[pos=0.4, above] {  $W_1^{(\varepsilon,\cof\varepsilon)}$};
\draw[->, thick] (8.2, -4) -- (8.2, -3.0); 

\draw[->, thick] (8.8, -4) -- (8.8, -3.0);

\draw[thick] (0.0, -4.5) -- (11.8, -4.5) node[pos=0.85, above] {  $W_{k-1}^{(\varepsilon,\cof\varepsilon)}$};
\draw[->, thick] (11.2, -4.5) -- (11.2, -0.8); 

\draw[->, thick] (11.8, -4.5) -- (11.8, -0.8);

\draw[thick] (0.0, -5) -- (13.8, -5) node[pos=0.905, above] {  $W_k^{(\varepsilon,\cof\varepsilon)}$};
\draw[->, thick] (13.2, -5) -- (13.2, -1.25); 

\draw[->, thick] (13.8, -5) -- (13.8, -1.25);
\end{tikzpicture}}
\caption{SPCNN architecture}
\end{figure}

Example activation functions $g_i$ are the ordinary maximum or the scaled $\mathrm{LogSumExp}$ functions 
$$\mathrm{LogSumExp}_r(x,y)=\frac{1}{r}\log(e^{rx}+e^{ry})$$ 
for $r>0$.
Note that as $r\to \infty$, $\mathrm{LogSumExp_r}\to \max$ and thus, the universal approximation property proved for $\max$-activation in Thm.~\ref{Thm: Universality of SPCNN in 3d} also extends to $\mathrm{LogSumExp_r}$ (for $r$ sufficiently large).
We also remark that choosing a grouping size of $2$ above is rather arbitrary, and was only done for simplicity. In theory, any number greater or equal than $2$ or even several distinct such grouping numbers as in \cite{tanielian2021approximating,anil2019sorting} could be chosen, we will, however, restrict ourselves to grouping size $2$ throughout this paper. 

\begin{Proposition}
    The lifted output function of the SPCNN is convex, and non-increasing in the second argument, i.e., the symmetric polyconvex output of the SPCNN is indeed symmetric polyconvex. 
\end{Proposition}
\begin{proof}
    Every $z_i$ in the iteration has the form 
    $$z_{i+1} =  g_i\biggl( W_i^{(z)}z_i +\bigl( W_{i,j}^{(\xi)}:\xi+W_{i,j}^{(\eta)}:\eta\bigr)_{j=1}^{2n_{i+1}}+b_i\biggr) .$$ The interior function is componentwise the sum of an affine function $W_{i,j}^{(\xi)}:\xi+W_{i,j}^{(\eta)}:\eta+b_{i,j}$ and a nonnegative linear combination of convex functions  $W_{i,j}^{(z)}z_i$, hence convex. Composing a convex function with a convex non-decreasing function produces a convex function, thus the lifted output is convex by induction. Concerning monotonicity: every component of $z_0=0$ is non-increasing in the second variable. Now, suppose that every component of $z_i$ is non-increasing in the second variable. Every nonnegative linear combination is necessarily non-increasing in the second variable. Moreover, the components of the affine function $\bigl( W_{i,j}^{(\xi)}:\xi+W_{i,j}^{(\eta)}:\eta\bigr)_{j=1}^{2n_{i+1}}+b_i$ are non-increasing in the second variable, therefore every component of the function $ W_i^{(z)}z_i +\bigl( W_{i,j}^{(\xi)}:\xi+W_{i,j}^{(\eta)}:\eta\bigr)_{j=1}^{2n_{i+1}}+b_i$ is non-increasing in the second argument. Composition with a non-decreasing function $g_i$ preserves this property. Thus the claim follows by induction over $i$.
\end{proof}

It is now straightforward to propose reduced architectures for the 2d case.

\begin{Definition}[SPCNN architecture in 2d]
        Let $\xi\in \Symtwo$ be an input matrix, let $k\in \N$, $n_1,...,n_k\in \N$, $n_0=n_{k+1}=1$, and let 
    \begin{align*}
        W_{i}^{(z)}&\in [0,\infty)^{2n_{i+1}\times n_{i}}, \quad b_i\in \R^{2n_{i+1}}\\
W_{i}^{(\xi)}&=\bigl(W_{i,j}^{(\xi)}\bigr)_{j=1}^{2n_{i+1}}\in \bigl(\Symtwo\bigr)^{2n_{i+1}} ,  \quad W_{i}^{(t)}=\bigl(W_{i,j}^{(t)}\bigr)_{j=1}^{2n_{i+1}}\in (-\infty,0]^{2n_{i+1}} 
    \end{align*}
be parameters for $i=0,...,k$, and $W_0^{(z)}=0$. We can group the $2n_{i+1}$ neurons into $n_{i+1}$ pairs of adjacent neurons. Furthermore, let 
$g_i\colon \R^2\to \R$ be nonlinear activation functions that are convex and non-decreasing applied to two adjacent components for $i=0,...,k$.
Then we iteratively set
\begin{align*}
    z_0=0, \ z_{i+1}& = g_i\biggl( W_i^{(z)}z_i +\bigl( W_{i,j}^{(\xi)}:\xi+W_{i,j}^{(t)}t\bigr)_{j=1}^{2n_{i+1}}+b_i\biggr), \  f(\xi,t,\Theta) = z_{k+1}.
\end{align*}Then the set $\Theta$ of all parameters together with the iteration scheme is called a 2d symmetric polyconvex neural network (\emph{SPCNN}) with activations $(g_i)_{i=0}^k$, $\mathcal{F}_{\mathrm{lift}}(\cdot,\Theta)\colon \Symtwo\times \R\to \R$ is called the lifted output function of the \emph{SPCNN} and $\mathcal{F}_{\mathrm{SPCNN}}(\cdot,\Theta)\colon \Symtwo\to \R, \varepsilon\mapsto \mathcal{F}_{\mathrm{lift}}(\varepsilon,\det\varepsilon,\Theta)$ is called the symmetric polyconvex output of the \emph{SPCNN}.
\end{Definition}

By the characterization of isotropic symmetric polyconvex functions, the representer function, defined on the diagonal matrices, admits certain symmetry properties. More precisely, the representer is invariant under permutation of the inputs. We can easily achieve this invariance for our network by simply averaging over all permutations. This does not interfere with the subsequent universal approximation property, cf.~\cite{geuken2025input}.
\begin{Definition}[SPCNN architecture for 3d isotropic functions]
       Let $x\in \R^3$ be an (eigenvalue) input, let $k\in \N$, $n_1,...,n_k\in \N$, $n_0=n_{k+1}=1$, and let 
    \begin{align*}
        W_{i}^{(z)}&\in [0,\infty)^{2n_{i+1}\times n_{i}}, \quad b_i\in \R^{2n_{i+1}}\\
W_{i}^{(x)}&=\bigl(W_{i,j}^{(x)}\bigr)_{j=1}^{2n_{i+1}}\in \bigl(\R^3\bigr)^{2n_{i+1}} ,  \quad W_{i}^{(y)}=\bigl(W_{i,j}^{(y)}\bigr)_{j=1}^{2n_{i+1}}\in \bigl((-\infty,0]^3\bigr)^{2n_{i+1}} 
    \end{align*}
be parameters for $i=0,...,k$, and $W_0^{(z)}=0$. We can group the $2n_{i+1}$ neurons into $n_{i+1}$ pairs of adjacent neurons. Furthermore, let 
$g_i\colon \R^2\to \R$ be nonlinear activation functions that are convex and non-decreasing applied to two adjacent components for $i=0,...,k$.
Then we iteratively set
\begin{align*}
    z_0=0, \ z_{i+1}& = g_i\biggl( W_i^{(z)}z_i +\bigl( W_{i,j}^{(x)}:x+W_{i,j}^{(y)}:y\bigr)_{j=1}^{2n_{i+1}}+b_i\biggr), \ f(x,y,\Theta) = z_{k+1}.
\end{align*}Then the set $\Theta$ of all parameters together with the iteration scheme is called a 3d isotropic symmetric polyconvex neural network (\emph{SPCNN}) with activations $(g_i)_{i=0}^k$, 
\begin{align*}
    \mathcal{F}_{\mathrm{lift}}(\cdot,\Theta)\colon \R^3\times \R^3&\to \R\\
    (x,y)&\mapsto \frac{1}{6}\biggl( \sum_{\sigma\in \mathrm{Perm}(3)}f(\sigma(x),\sigma(y),\Theta)\biggr)
\end{align*}
 is called the lifted output function of the \emph{SPCNN} and $\mathcal{F}_{\mathrm{SPCNN}}(\cdot,\Theta)\colon \R^3\to \R, x\mapsto \mathcal{F}_{\mathrm{lift}}(x,\cof x,\Theta)$ is called the symmetric polyconvex output of the \emph{SPCNN}.
\end{Definition}

\begin{Definition}[SPCNN architecture for 2d isotropic functions]
            Let $x\in \R^2$ be an (eigenvalue) input, let $k\in \N$, $n_1,...,n_k\in \N$, $n_0=n_{k+1}=1$, and let 
    \begin{align*}
        W_{i}^{(z)}&\in [0,\infty)^{2n_{i+1}\times n_{i}}, \quad b_i\in \R^{2n_{i+1}}\\
W_{i}^{(x)}&=\bigl(W_{i,j}^{(x)}\bigr)_{j=1}^{2n_{i+1}}\in \bigl(\R^2\bigr)^{2n_{i+1}} ,  \quad W_{i}^{(t)}=\bigl(W_{i,j}^{(t)}\bigr)_{j=1}^{2n_{i+1}}\in (-\infty,0]^{2n_{i+1}} 
    \end{align*}
be parameters for $i=0,...,k$, and $W_0^{(z)}=0$. We can group the $2n_{i+1}$ neurons into $n_{i+1}$ pairs of adjacent neurons. Furthermore, let 
$g_i\colon \R^2\to \R$ be nonlinear activation functions that are convex and non-decreasing applied to two adjacent components for $i=0,...,k$.
Then we iteratively set
\begin{align*}
    z_0=0, \ z_{i+1}& = g_i\biggl( W_i^{(z)}z_i +\bigl( W_{i,j}^{(x)}:x+W_{i,j}^{(t)}t\bigr)_{j=1}^{2n_{i+1}}+b_i\biggr), \ f(x,t,\Theta) = z_{k+1}.
\end{align*}Then the set $\Theta$ of all parameters together with the iteration scheme is called a 2d isotropic symmetric polyconvex neural network (\emph{SPCNN}) with activations $(g_i)_{i=0}^k$, 
\begin{align*}
    \mathcal{F}_{\mathrm{lift}}(\cdot,\Theta)\colon \R^2\times \R&\to \R\\
    (\lambda_1,\lambda_2,t)&\mapsto \frac{1}{2}\bigl( f(\lambda_1,\lambda_2,t,\Theta)+f(\lambda_2,\lambda_1,t,\Theta)\bigr)
\end{align*}
 is called the lifted output function of the \emph{SPCNN} and $\mathcal{F}_{\mathrm{SPCNN}}(\cdot,\Theta)\colon \R^2\to \R, (\lambda_1,\lambda_2)\mapsto \mathcal{F}_{\mathrm{lift}}(\lambda_1,\lambda_2,\lambda_1\lambda_2,\Theta)$ is called the symmetric polyconvex output of the \emph{SPCNN}.
\end{Definition}

\subsection{Universality and expressivity of SPCNN}

In the following part of the section, we prove that our 3d architecture has the universal approximation property. The proofs for the 2d and isotropic cases are analogous and therefore omitted. The procedure is inspired by \cite{huang2021convex,geuken2025input}.

\begin{lemma}\label{Approximation of symmetric polyconvex functions by maximum of affine functions}
    Let $U,V\subset \Symthree$ be nonempty and compact. Let $g\colon \Symthree\times \Symthree\to \R$ be convex and non-increasing in the second variable. For any $\tau>0$ there are $k$ affine functions $L_j\colon \Symthree\times \Symthree\to \R$, $j=1,\ldots, k$, $k\in \N$, of the form
    $$L_j(\xi,\eta)=\alpha_j:\xi+\beta_j:\eta+c_j,$$
    where $\alpha_j\in \Symthree, \beta_j\in \Snegativthree, c_j\in \R$, such that
    $$\sup_{(\xi,\eta)\in U\times V}\biggl\vert g(\xi,\eta)-\max_{j=1,\ldots,k}(L_j)(\xi,\eta)\biggl\vert\le \tau.$$
\end{lemma}
\begin{proof}Let $\tau>0$. Let us equip the product space with the norm $\Vert (\xi,\eta)\Vert:=\sqrt{\vert \xi\vert^2+\vert\eta\vert^2}$.
    Since $g$ is finite and convex, it is continuous, its subdifferential is nonempty everywhere, and its subgradients are locally bounded. By compactness of $U\times V$,
    $$M:=\sup\biggl\{ \Vert (\alpha,\beta)\Vert : (\xi,\eta)\in U\times V,\ (\alpha,\beta)\in\partial g(\xi,\eta)\biggr\}<\infty.$$
Uniform continuity of $g$ allows us to choose $\delta>0$ such that $M\delta<\frac{\tau}{2}$ and
$$\vert g(\xi_1,\eta_1)-g(\xi_2,\eta_2)\vert <\frac{\tau}{2}
\quad \text{whenever }\xi_1,\xi_2\in U,\eta_1,\eta_2\in V,\ \Vert (\xi_1,\eta_1)-(\xi_2,\eta_2)\Vert<\delta.$$
Choose $\xi_j\in U,\eta_j \in V$, $j=1,\dots,k$, such that
$$U\times V\subset\bigcup_{j=1}^k B_\delta(\xi_j,\eta_j),$$
and choose
$$(\alpha_j,\beta_j)\in\partial g(\xi_j,\eta_j).$$
Define the supporting affine functions
$$L_j(\xi,\eta) =g(\xi_j,\eta_j) +\alpha_j:(\xi-\xi_j) +\beta_j:(\eta-\eta_j).$$
By the subgradient inequality, $L_j\le g$ globally.
Note that by monotonicity $\beta_j\in \Snegativthree$.
Finally, given $\xi\in U$, $\eta\in V$, choose $j$ with $\Vert(\xi,\eta)- (\xi_j,\eta_j)\Vert<\delta$. Then
\begin{align*}
0\le g(\xi,\eta)-\max_{1\le i\le k}L_i(\xi,\eta) &\le g(\xi,\eta)-L_j(\xi,\eta)\\ &=g(\xi,\eta)-g(\xi_j,\eta_j)-\alpha_j:(\xi-\xi_j)-\beta_j:(\eta-\eta_j)\\ 
&\le \vert g(\xi,\eta)-g(\xi_j,\eta_j)\vert +\Vert (\alpha_j,\beta_j)\Vert \Vert (\xi,\eta)-(\xi_j,\eta_j)\Vert \\ 
&<\tau.
\end{align*}
\end{proof}

\begin{lemma} \label{Approximate lifted polyconvex fucntions}
 Let $U,V\subset \Symthree$ be nonempty and compact. Let $g\colon \Symthree\times \Symthree\to \R$ be convex and non-increasing in the second variable. For any $\tau>0$ there is a \emph{SPCNN} with $\max$-activation $\mathcal{F}(\cdot,\Theta)$ such that
 $$\sup_{(\xi,\eta)\in U\times V}\vert g(\xi,\eta)-\mathcal{F}_{\mathrm{lift}}(\xi,\eta,\Theta)\vert\le \tau.$$
\end{lemma}
\begin{proof}
By lemma \ref{Approximation of symmetric polyconvex functions by maximum of affine functions}
    it suffices to show the claim for $g$ of the form
    $$ g(\xi,\eta)=\max_{j=0,\ldots,k}(L_j)(\xi,\eta)$$
    where $k\in \N$ and $L_j(\xi,\eta)=\alpha_j:\xi+\beta_j:\eta+c_j$,
     $\alpha_j\in \Symthree, \beta_j\in \Snegativthree, c_j\in \R$. Such piecewise affine convex functions can be written 
     as 
     \begin{align*}
         g(\xi,\eta)=\max(L_k,\max(L_{k-1},\max(...,\max(L_1,L_0)))).
     \end{align*}For the network, we may choose 
    $W_{j,1}^{(\xi)}=W_{j,1}^{(\eta)}=0, b_{j,1}=0$ and
     $W_{j,2}^{(\xi)},W_{j,2}^{(\eta)},b_{j,2}$ such that 
$$L_{j+1}(\xi,\eta)=W_{j,2}^{(\xi)}:\xi+W_{j,2}^{(\eta)}:\eta+b_{j,2}$$
     for all $\xi,\eta\in \Symthree$, $j=1,...,k-1$ 
  and choose $W_{0,i+1}^{(\xi)},W_{0,i+1}^{(\eta)},b_{0,i+1}$ such that
$$L_i(\xi,\eta)=W_{0,i+1}^{(\xi)}:\xi+W_{0,i+1}^{(\eta)}:\eta+b_{0,i+1}$$
for $i=0,1$.
     We set $W_{j,1}^{(z)}=1$ and $W_{j,2}^{(z)}=0$ for $j=1,...,k-1$. Then
     $$z_1=\max(L_1,L_0)(\xi,\eta).$$
     Furthermore,
     \begin{align*}
         z_{j+1}=\max(L_{j+1},z_{j}) .
     \end{align*}In particular, $z_{k+1}=g$.
\end{proof}

\begin{Theorem}[Universality of SPCNN in 3d]\label{Thm: Universality of SPCNN in 3d}
    Let $F\colon \Symthree\to \R$ be symmetric polyconvex, $U\subset \Symthree$ compact, and $\tau>0$. Then there is a 3d \emph{SPCNN} $\mathcal{F}(\cdot,\Theta)$ with $\max$-activation, such that
    $$\sup_{\varepsilon\in U}\vert F(\varepsilon)-\mathcal{F}_{\mathrm{SPCNN}}(\varepsilon,\Theta)\vert\le \tau.$$
\end{Theorem}
\begin{proof}
    By the characterization of symmetric polyconvex functions in Thm.~\ref{Kreisbeck Theorem}, there is a representer function
    $$G\colon \Symthree\times \Symthree\to \R$$
    that is convex, non-increasing in the second argument and $G(\varepsilon,\cof\varepsilon)=F(\varepsilon)$ for all $\varepsilon\in \Symthree$. For $\tau>0$, we can approximate $G$ by an SPCNN $f(\cdot, \Theta)$ according to lemma \ref{Approximate lifted polyconvex fucntions}. Precisely, let $V=\cof(U)$, which is compact by the continuity of the cofactor mapping. Then
    \begin{align*}
        \sup_{(\xi,\eta)\in U\times V}\vert G(\xi,\eta)-\mathcal{F}_{\mathrm{lift}}(\xi,\eta,\Theta)\vert\le \tau,
    \end{align*}
    and by construction
    \begin{align*}
        \vert F(\varepsilon)-\mathcal{F}_{\mathrm{SPCNN}}(\varepsilon,\Theta)\vert &= \vert G(\varepsilon,\cof\varepsilon)-\mathcal{F}_{\mathrm{lift}}(\varepsilon,\cof\varepsilon,\Theta)\vert \le \tau.
    \end{align*}
\end{proof}

We next quantify the expressivity of our SPCNN architecture by comparing the number of its parameters with the number of affine pieces required in an equivalent single-layer max-out representation.

\begin{Theorem}[Expressivity, cf.~\cite{chen2019optimalcontrolneuralnetworks}, Thm.~\ref{Thm: Expressivity of ICNN}]
   There is a single hidden layer \emph{SPCNN} with $\max$-activation and $2N$ neurons, $N\ge 12$, such that there are pairwise distinct affine functions $L_i\colon \Symthree\times \Symthree\to \R$, $i=1,...,C(N,12)$ that are non-increasing in the second argument, such that 
$$\mathcal{F}_{\mathrm{lift}}(\cdot, \Theta)=\max_{i=1,...,C(N,12)}L_i,$$
   where
 $$C(N,12):=\sum_{j=0}^{12}\binom{N}{j}=\frac{N^{12}}{12!}+\mathcal{O} (N^{11}).$$
\end{Theorem}

\begin{proof}
    Take an SPCNN with 1 hidden layer and $2N$ neurons and choose the affine functions generically like in Thm.~\ref{Thm: Expressivity of ICNN}. Set $W_{i,j}^{(z)}=1$ for $j=1,...,N$. Then
    $$z_1=(\max(L'_{2j-1},L'_{2j}))_{j=1}^N$$
    for the affine functions resulting from the weights and biases of the hidden neurons. 
    We set $W_{1,1}^{(\xi)}=W_{1,1}^{(\eta)}=b_{1,1}=0$, and we proceed similarly for the second component of $z_2$. Then
    $$\mathcal{F}_{\mathrm{lift}}(\cdot,\Theta)=\max(z_{2,1},z_{2,2})=z_{2,1}=\sum_{j=1}^N \max(L'_{2j-1},L'_{2j}).$$
    As in Thm.~\ref{Thm: Expressivity of ICNN}, the affine functions can be made pairwise distinct, rendering a separation of the 12 dimensional space $\Symthree\times \Symthree$ using the $N$ hyperplanes $\{L'_{2j-1}=L'_{2j}\}$. On each separated region, $\mathcal{F}_{\mathrm{lift}}(\cdot,\Theta)$ is affine and different from any other region, satisfying the desired property. 
\end{proof}
\begin{Bem}
    The construction in the above proof also shows, by setting every even $L_j$ to zero, that our \emph{SPCNN}s contain the regular \emph{ICNN}s (with restricted parameters in the second argument) as a special case. 
\end{Bem}
Next, approximating $F\circ \pi_S$ by a polyconvex network requires considering $9$ instead of $6$ inputs. This approach leads to a larger number of parameters, and thus more compute, compared to our SPCNNs, which make use of the direct characterization of symmetric polyconvex functions in terms of symmetric matrices. 

Finally, one could use a standard ICNN to approximate a symmetric convex representer function $\tilde F$ of $F$ and enforce the monotonicity either with a penalization in the cost-function, preconditioning the training data, or directly by sorting out violations, e.g., if the network produces any polyaffine parts $\beta:\cof\varepsilon$ for $\beta\notin \Snegativthree$, then we project $\beta$ onto the cone of symmetric negative semidefinite matrices. 

The first approach neither guarantees monotonicity nor does it require fewer parameters than our SPCNN. 
The second approach uses training pairs (see Section 4.1) filtered to satisfy the required monotonicity in the cofactor arguments, simply by discarding any non-monotonic pairs from the dataset.
This, again, does not yield monotonicity of the output function, but only slashes the training set and could lead to a worse approximation of the symmetric polyconvex envelope in the end. 
The main point of critique of the third approach is that it potentially produces constant parts in the output that can be far from the desired outcome.

\section{Applications of SPCNNs}

In this section, we present two possible applications of our symmetric polyconvex neural networks for engineering purposes. As a general remark, the SPCNN approach has the advantage, compared to standard numerical methods, that the output is symmetric polyconvex, and thus the associated energy minimizing problem admits a solution (under additional suitable growth conditions).

\subsection{Identification of hyperelastic energies for experimental data}

Quadratic elastic energies describe a linear stress–strain response, whereas nonlinear elastic behavior requires more flexible constitutive representations. 
Our SPCNNs enable the formulation of such nonlinear material models through a statistical approach. We assume that only a limited number of experimental stress-strain pairs are available for a given material. Generalizing from this sparse training set to make predictions outside the training set poses a classical regression problem, which we address using SPCNNs. Moreover, our SPCNNs produce a global closed analytical formula, which is different from a piecewise interpolation. 

We note that the universal approximation theorem \ref{Thm: Universality of SPCNN in 3d} extends to simultaneous approximation of $C^1$-energies and stresses on compact subsets by means of the following lemma.
\begin{lemma}
Let $U\subset\Symthree$ be open, let $F\in C^1(U,\R)$, and let $F_j\in C^1(U,\R)$ be convex along every rank-one direction $a\otimes a$, $a\in \R^3$. If
$$F_j\longrightarrow F
\quad \text{locally uniformly in } U,$$
then
$$\nabla F_j\longrightarrow\nabla F
\quad \text{locally uniformly in } U.$$
In particular, this applies to smooth symmetric polyconvex functions.
\end{lemma}
\begin{proof}
    We adapt the difference-quotient argument in the proof of \cite[Thm~25.7]{Rockafellar97convex}. For this, convexity is only needed along the direction being tested, and finitely many $a\otimes a$ span $\Symthree$.
\end{proof}

We will take synthetic data here by sampling $n$ stress-strain pairs from a hyperelastic energy density consisting of a classical quadratic elastic and a quartic hardening term:
$$F(\varepsilon)=\frac{\kappa}{2}(\tr\varepsilon)^2+\mu\Vert\dev\varepsilon\Vert^2+\alpha \Vert\dev\varepsilon\Vert^4.$$
Since in applications we do not have $F$ at hand, but rather its gradient, the stress, our training set consists of strain-stress pairs
$$\mathscr{T}_1=\{(\varepsilon_i,\sigma(\varepsilon_i))\mid i=1,...,n\},$$
where
$$\sigma(\varepsilon)=\nabla F(\varepsilon)=\kappa  \tr\varepsilon\mathbb{I}+2\mu\dev  \varepsilon+4\alpha \Vert \dev \varepsilon\Vert^2\dev \varepsilon.$$
Next, we need to minimize the distance between the derivative of $F$ (the measured stresses) and the derivative of the SPCNN. The gradient of our SPCNN 
$$\mathcal{P}_{\mathrm{SPCNN}}(\varepsilon,\Theta):=\frac{\partial}{\partial \varepsilon} \mathcal{F}_{\mathrm{SPCNN}}(\varepsilon,\Theta)$$
can efficiently be calculated using backpropagation. 
Now, it remains to solve the Sobolev regression problem. 
As a cost function, we take the squared error
$$\mathrm{Cost}_2 (x,y)=\frac{1}{2}\vert x-y\vert^2 $$
The loss function for our problem is then the mean squared error
$$L_1(\Theta)=\frac{1}{2n}\sum_{i=1}^n\bigl\vert \mathcal{P}_{\mathrm{SPCNN}}(\varepsilon_i,\Theta)-\sigma(\varepsilon_i)\bigl\vert^2.$$
We choose the analogous loss function for the 2d and the isotropic cases. 
We will always use the SPCNN architecture with pairwise $\mathrm{LogSumExp}_r$-activation (with $r=20$), except for the last layer which has a linear output. Note that the result of Sobolev training is not unique, therefore we employ an additional normalization $\mathcal{F}_{\mathrm{SPCNN}}(0,\Theta)=0$.

\subsection{SPCNN for computing symmetric polyconvexification}

After explaining how our SPCNNs can be used to represent experimental data in an energetic framework, it is now in order to demonstrate how our SPCNNs can be used to approximate the symmetric polyconvex envelope 
$$F^{\mathrm{spc}}(\varepsilon)=\sup\{G(\varepsilon)\mid \ G \ \text{is symmetric polyconvex, }G\le F\}$$
of a given function $F\colon \Symthree\to \R$, which may be derived from microscopic considerations. To approximate $F^{\mathrm{spc}}$ with an SPCNN, we first need to choose an appropriate loss function, such that the resulting symmetric polyconvex output lies as close to $F$ as possible (regression problem) while simultaneously penalizing violations of the upper bound assumption. To achieve this, we proceed along the lines of \cite[Sec.~5.3]{geuken2025input}. To make stochastic gradient descent possible, we choose a loss function that is an average over the individual costs per training datum. For the regression part, there are multiple reasonable cost functions, such as the absolute value of the difference, the square of the difference, and so on. Crucially, we remark that it is not the purpose of the envelope to solve the regression problem exactly. More precisely, the difference between the function and its envelope can be very large in the $L^{\infty}$-norm. This means that taking the squared distance gives greater weight to large deviations and may therefore require a different balance between the fitting and penalty terms, hence we opt for the absolute value in this paper. For the penalization, we naturally choose $\relu$ of the difference of network output and given function $F$. To summarize, our cost function is
$$\mathrm{Cost}_{1,s} (x,y)=\vert x-y\vert +s \relu(x-y),$$
where $s$ is an a priori chosen hyperparameter. 
Our training set is a sample of the graph of $F$, meaning, we take $n$ distinct tensors $\varepsilon_i \in \Symthree$, $i=1,...,n$, and the training set is
$$\mathscr{T}_2=\{(\varepsilon_i,F(\varepsilon_i))\mid i=1,...,n\}.$$
The loss function for our problem is then
$$L_{1,s}(\Theta)=\frac{1}{n}\sum_{i=1}^n \biggl[\bigl\vert \mathcal{F}_{\mathrm{SPCNN}}(\varepsilon_i,\Theta)-F(\varepsilon_i)\bigl\vert + s \relu\bigl(\mathcal{F}_{\mathrm{SPCNN}}(\varepsilon_i,\Theta)-F(\varepsilon_i)\bigr)\biggr].$$
We choose the analogous loss function for the 2d and the isotropic cases. 

The example we consider is a classical double-well potential. We assume that our material exhibits two phases, each of which admits the elastic energy
$$\frac{\kappa_i}{2}(\tr\varepsilon)^2+s_i \Vert \dev \varepsilon\Vert^2+\gamma_i.$$
We will use the SPCNN architecture with pairwise $\mathrm{LogSumExp}_r$-activation (with $r=20$), except for the last layer which has a linear output. 
For the isotropic case, we will set the last pass-through layer to zero.

\section{Numerical results}

Regarding the general training procedure, the most fundamental gradient-based approach for solving the weight optimization problem is projected subgradient descent, alongside variants that incorporate a momentum term.

All models are implemented in Python using the PyTorch framework, making use of the Adam optimizer \cite{kingma2015adam}.

\subsection{Sobolev training}

We consider a two-dimensional isotropic example with $\kappa=2.0, \mu=6.0, \alpha=4.0$
$$f(x,y)=(x+y)^2+3(x-y)^2+(x-y)^4.$$
Then 
\begin{align*}
    \sigma_x(x,y)&=2(x+y)+6(x-y)+4(x-y)^3\\
     \sigma_y(x,y)&=2(x+y)-6(x-y)-4(x-y)^3.
\end{align*}
We use an equi-distant $10\times 10$-grid in $[-2,2]^2$. 
Our network architecture has three hidden layers of width (before activation) $12, 18, 6$ respectively. We train, using the PyTorch Adam implementation, for 1000 epochs with a learning rate of 0.005. The training lasts for roughly 2 seconds on a standard laptop. 
The SPCNN accurately captures the complex quartic growth in the deviatoric direction. While the approximation of the quadratic trace growth appears less precise at first glance, this is a direct consequence of the orders-of-magnitude difference in scale between the two directions in this benchmark. Accounting for these scale disparities, the overall approximation is robust. 
\begin{figure}[H]
    \centering
        \includegraphics[scale=0.6]{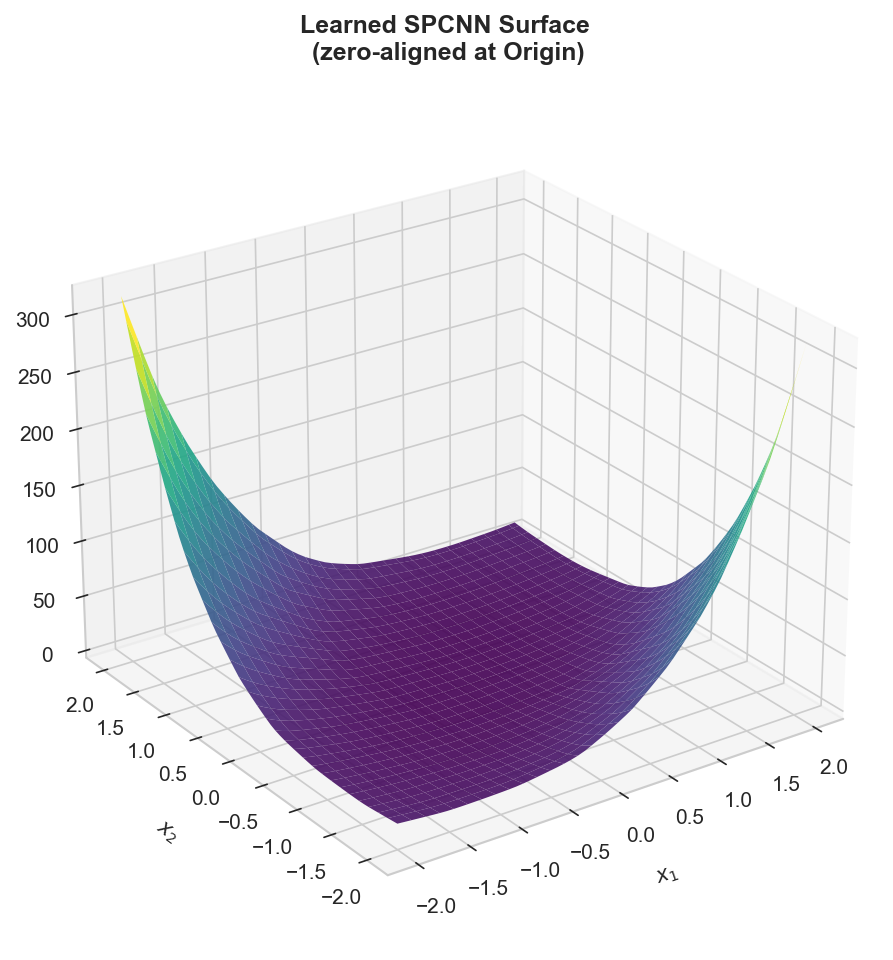}
    \caption{SPCNN output surface}
    \label{Fig: SPCNN output Surface}
\end{figure}

\begin{figure}[H]
        \centering
        \includegraphics[scale=0.4]{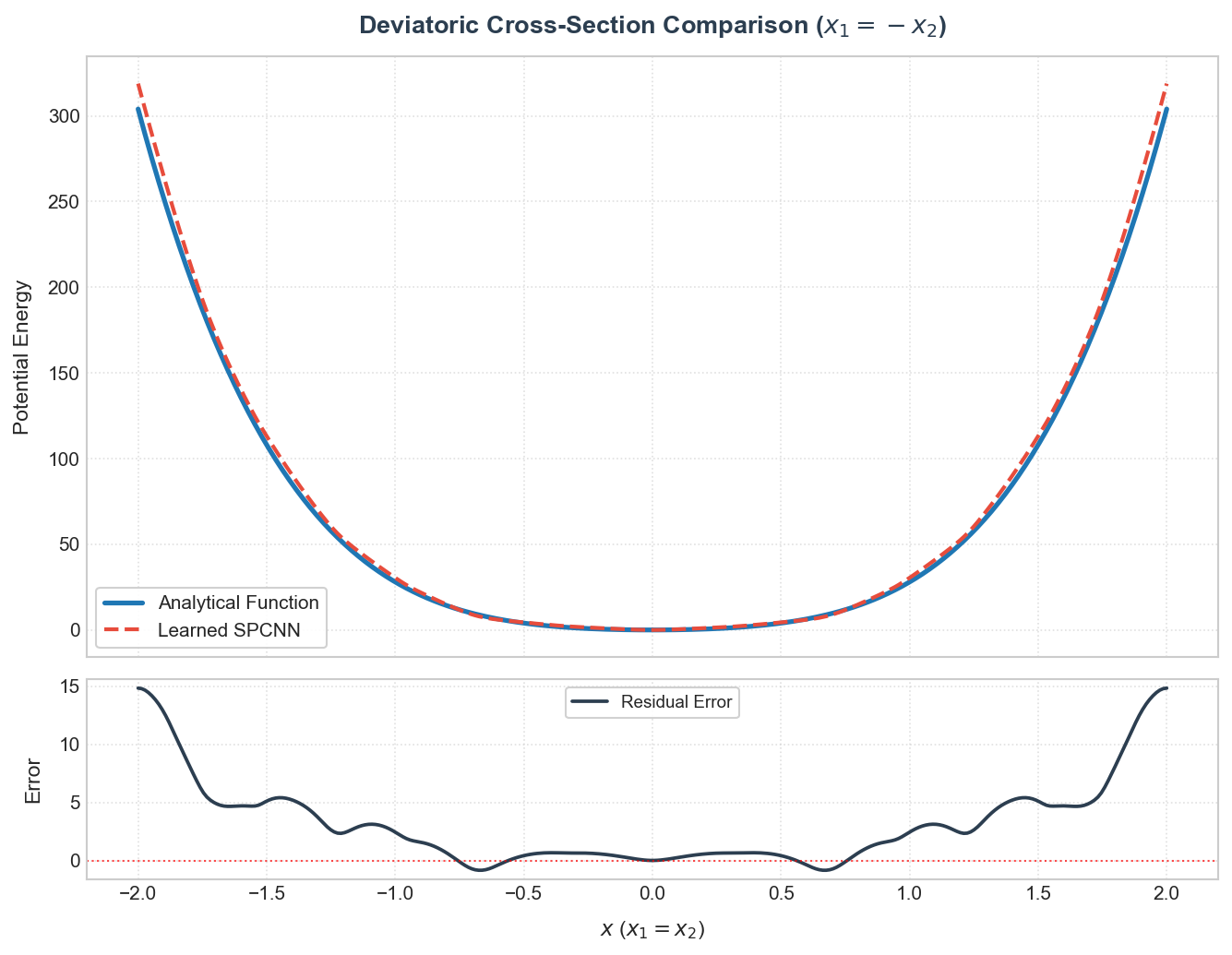}
    \caption{SPCNN output on the deviatoric line}
    \label{Fig: Sobolev Training deviatoric}
\end{figure}
 While fidelity can be systematically improved by scaling the network capacity, providing more training data, or tuning hyperparameters, our baseline configuration shows that even a small amount of training yields effective and practically useful results.
\begin{figure}[H]
        \centering
        \includegraphics[scale=0.4]{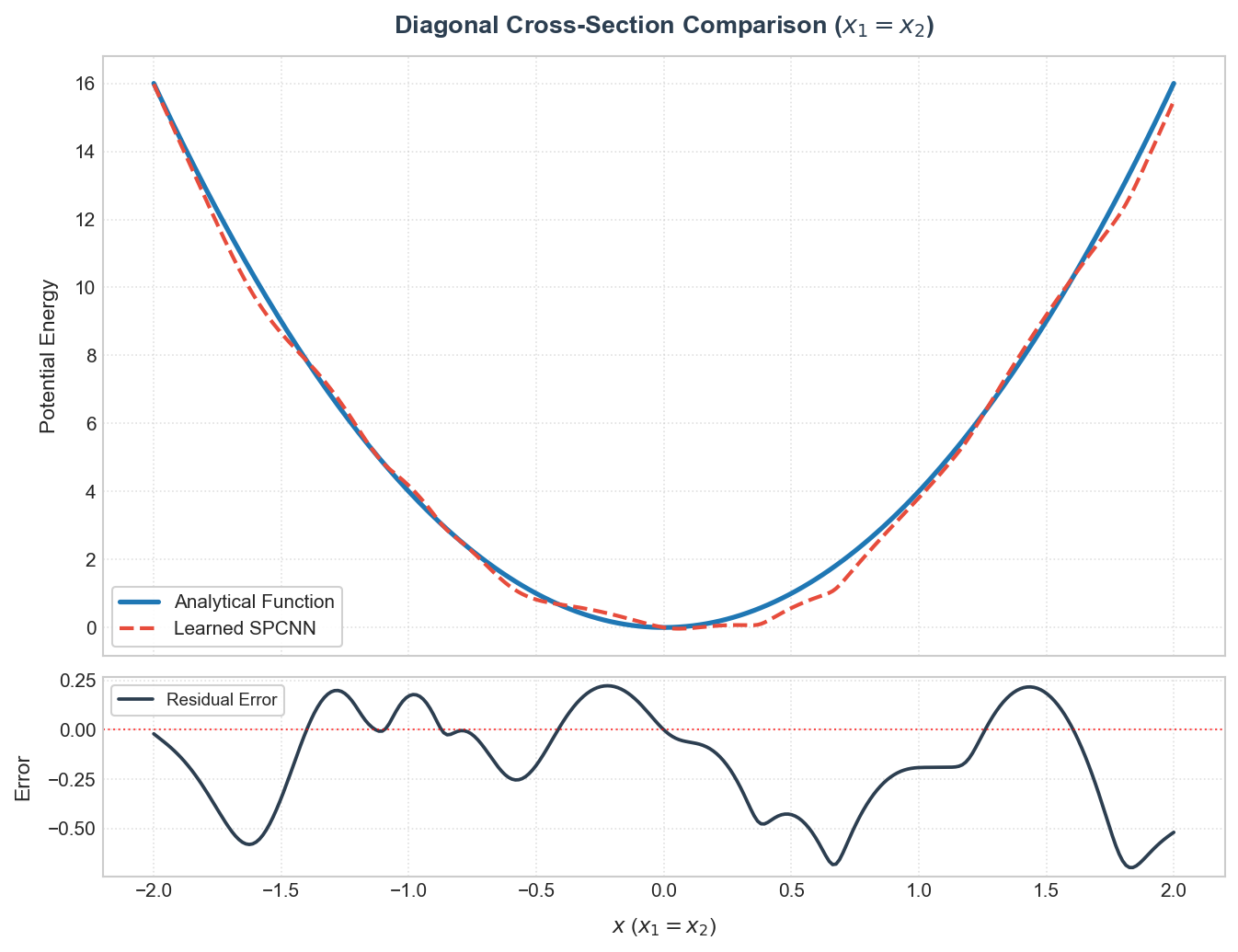}
    \caption{SPCNN output on the diagonal line}
    \label{Fig: Sobolev Training diagonal}
\end{figure}

\subsection{Relaxation of a double-well type problem}

The example we consider here is an isotropic two-dimensional double-well potential, where we deliberately chose the normalized model parameters to differ greatly in order to produce substantial differences between the energy and its symmetric polyconvex and convex envelopes. 
We choose 
\begin{align*}
    f_1(x,y)&= 40(x+y)^2+18(x-y)^2, \\
    f_2(x,y)&= (x+y)^2 +(x-y)^2 +40,\\
    f(x,y)&:=\min(f_1,f_2)(x,y)
\end{align*}
For this type of function, the symmetric polyconvex  envelope is known explicitly (see \cite{DolGam2026}). It differs the most from the convex envelope along $x=y$ (the trace/hydrostatic line). 
Along this line, the envelope reads for $x\ge 0$
\begin{align*}
    f^{\mathrm{spc}}(x,x)&= \begin{cases} f_1(x,x) &\text{if }2x\le \sqrt{\frac{160}{82\cdot 39}},\\
   f_2(x,x) &\text{if }2x\ge \sqrt{\frac{820}{ 39}},\\
        164\sqrt{\frac{160}{82\cdot 39}}x-\frac{80}{39} -4x^2 & \text{otherwise}.
    \end{cases}
\end{align*}
For the training, we use a grid that is not linear but biased towards the regions where the desired function, $f^{\mathrm{spc}}$, coincides with $f_1$ (small values of $x,y$), and with $f_2$ (large values of $x,y$) respectively. 
We do this by first creating a set of 100 equi-distant points in $[-2,2]$, then raise every such point to the power of $5$, do this for the $x$- and $y$-direction, and then form the Cartesian product to get an adapted mesh of 10,000 points on $[-32,32]^2$. Our penalization parameter is set to $s=r=20$.

\begin{figure}[H]
        \centering
        \includegraphics[scale=0.52]{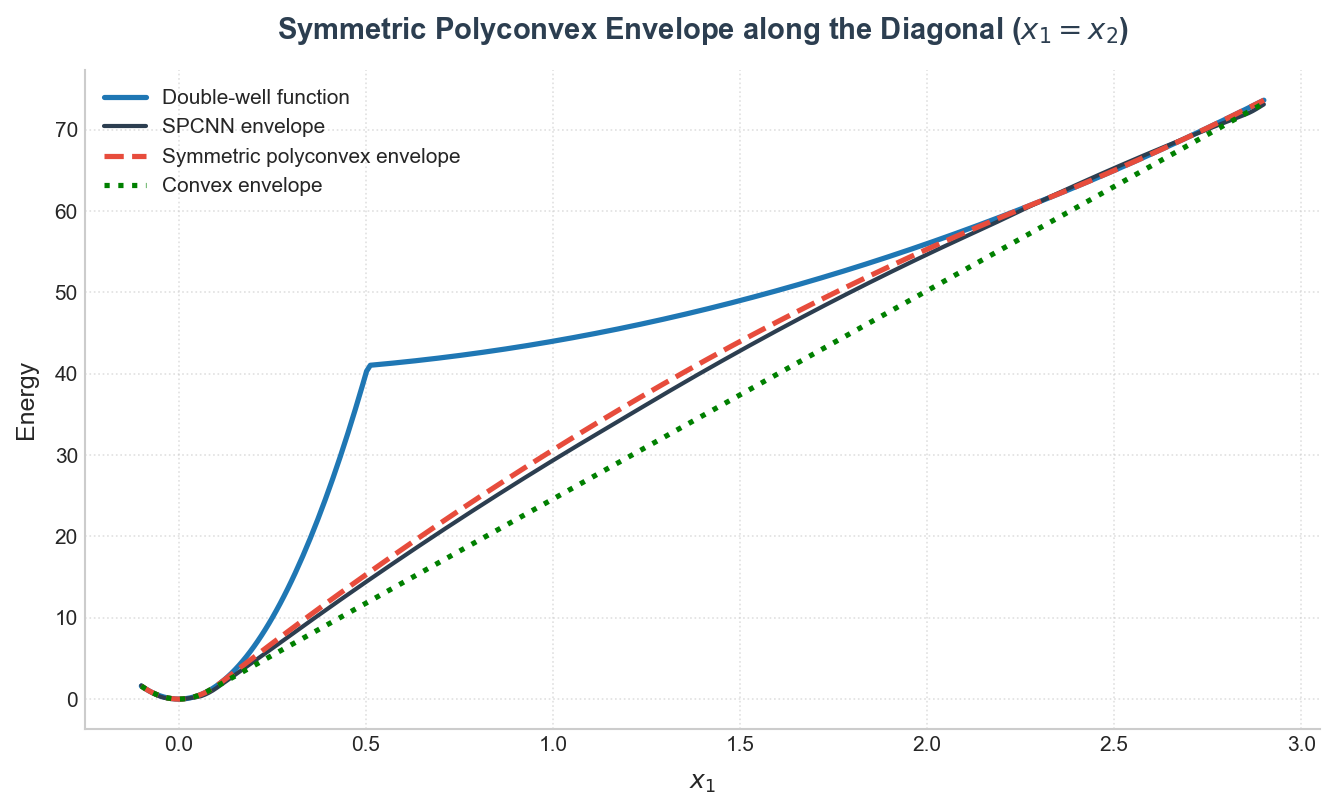}
    \caption{SPCNN envelope along the diagonal}
    \label{Fig: Envelope-Training Output Comparison}
\end{figure}
Our network architecture has four hidden layers of width (before activation) $12, 18, 12, 6$. 
The training lasts for 1500 epochs with a constant batch size of 256, employing the PyTorch Adam optimizer. We employ a piecewise learning rate schedule, updating the rate every 300 epochs as follows:

\begin{enumerate}[noitemsep]
    \item[1.] 300 epochs with a learning rate of 0.005,
    \item[2.] 300 epochs with a learning rate of 0.001,
    \item[3.] 300 epochs with a learning rate of 0.0005,
    \item[4.] 300 epochs with a learning rate of 0.0001,
    \item[5.] 300 epochs with a learning rate of 0.00001.
\end{enumerate}
All numerical experiments were conducted on a standard laptop and required approximately three minutes of compute time.

We see from our result that the approximation for smaller strains is highly accurate. The relative error is within a 5\% margin on the whole interval considered (apart from a small vicinity of the origin caused by a vanishing denominator), and our SPCNN output is a significantly better approximation than the convex envelope.

At larger strains, the accuracy of the approximation starts to deteriorate.
This discrepancy arises because the neural network is trained on a bounded domain, whereas the analytical method considers the entire unbounded mathematical domain $\R^2$ (in general $\Sym$).
 Consequently, the network's predictions are tied to its training limits, while the analytical envelope remains completely domain-independent.

\begin{figure}[H]
        \centering
        \includegraphics[scale=0.6]{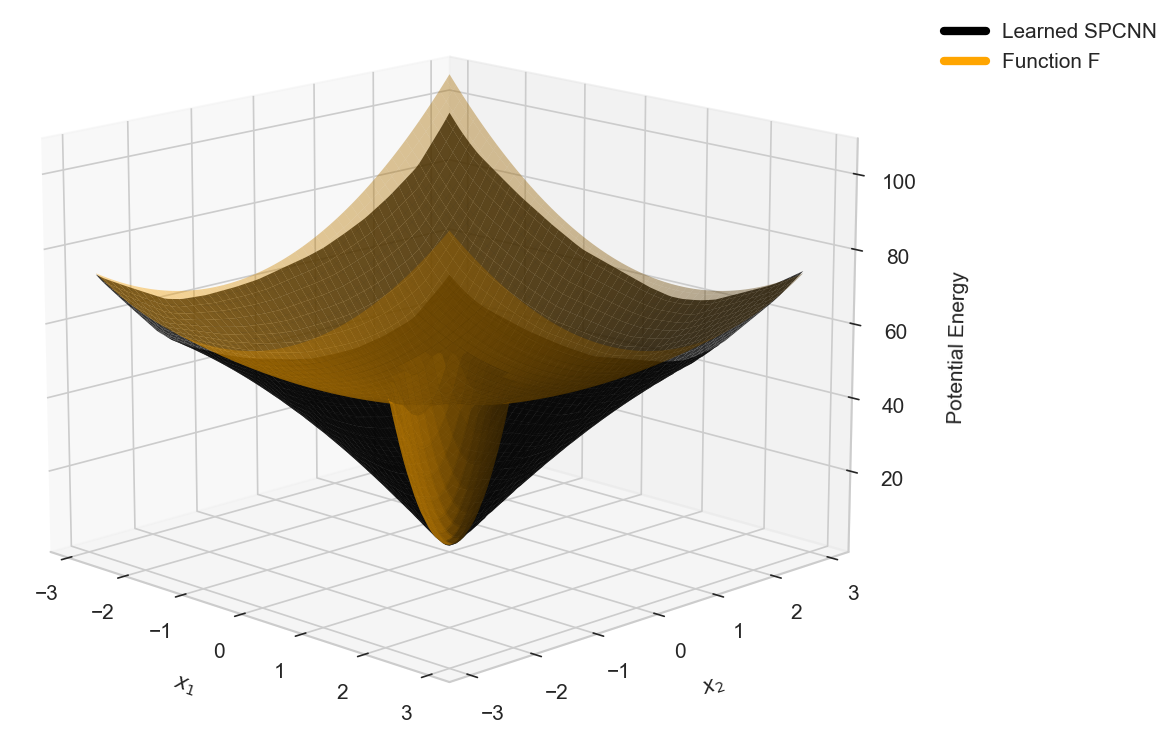}
    \caption{SPCNN envelope}
    \label{Fig: total Envelope-Training Output Comparison}
\end{figure}

However, this is not a particular issue of our neural network approach, but a general problem of numerically determining polyconvex envelopes. 

 Finally, the second condition of lying below $F$ is satisfied up to a relative error of less than 1\% for $x\in(2.3,3.6)$. Outside of this interval, the obstacle constraint is satisfied. 

\begin{figure}[H]
    \centering
    \begin{minipage}{0.5\textwidth}
        \centering
        \includegraphics[width=\textwidth]{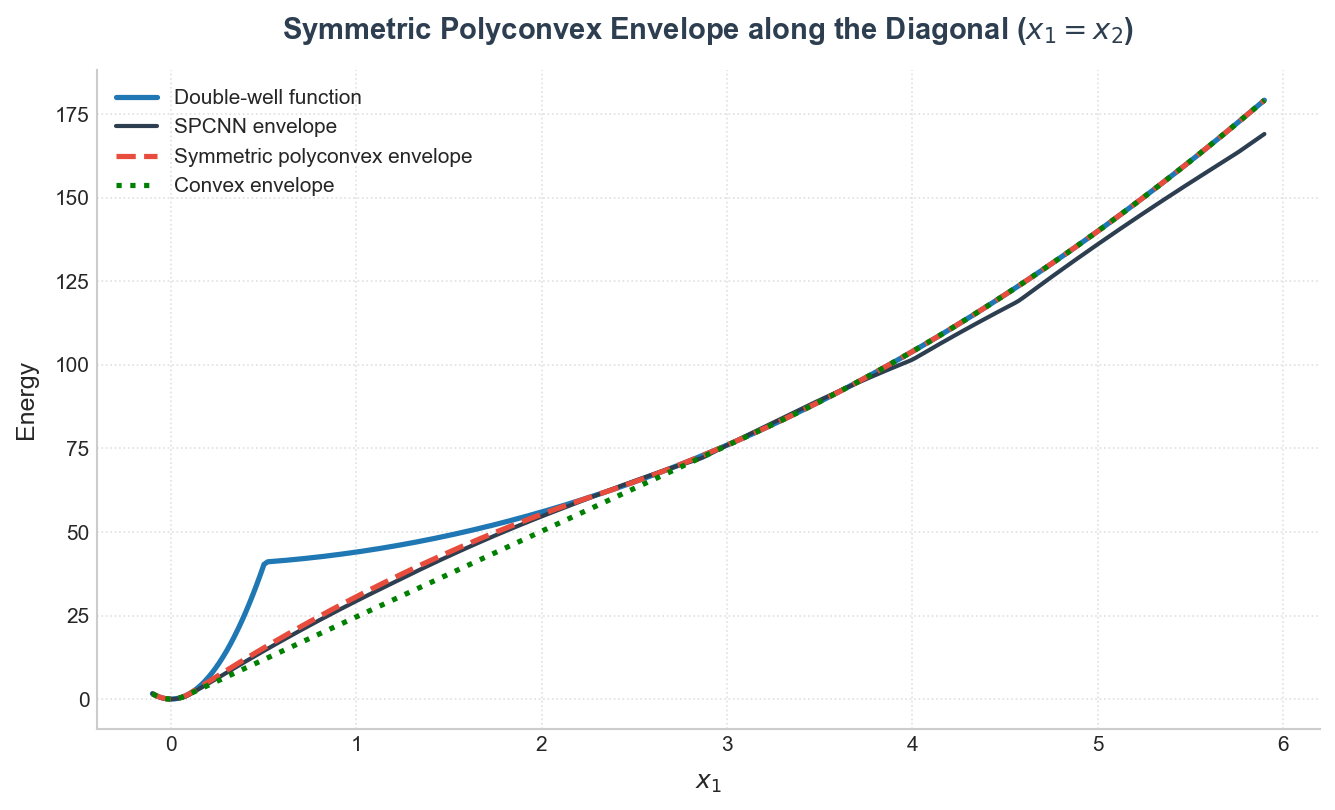}
    \end{minipage}\hfill
    \begin{minipage}{0.5\textwidth}
        \centering
        \includegraphics[width=\textwidth]{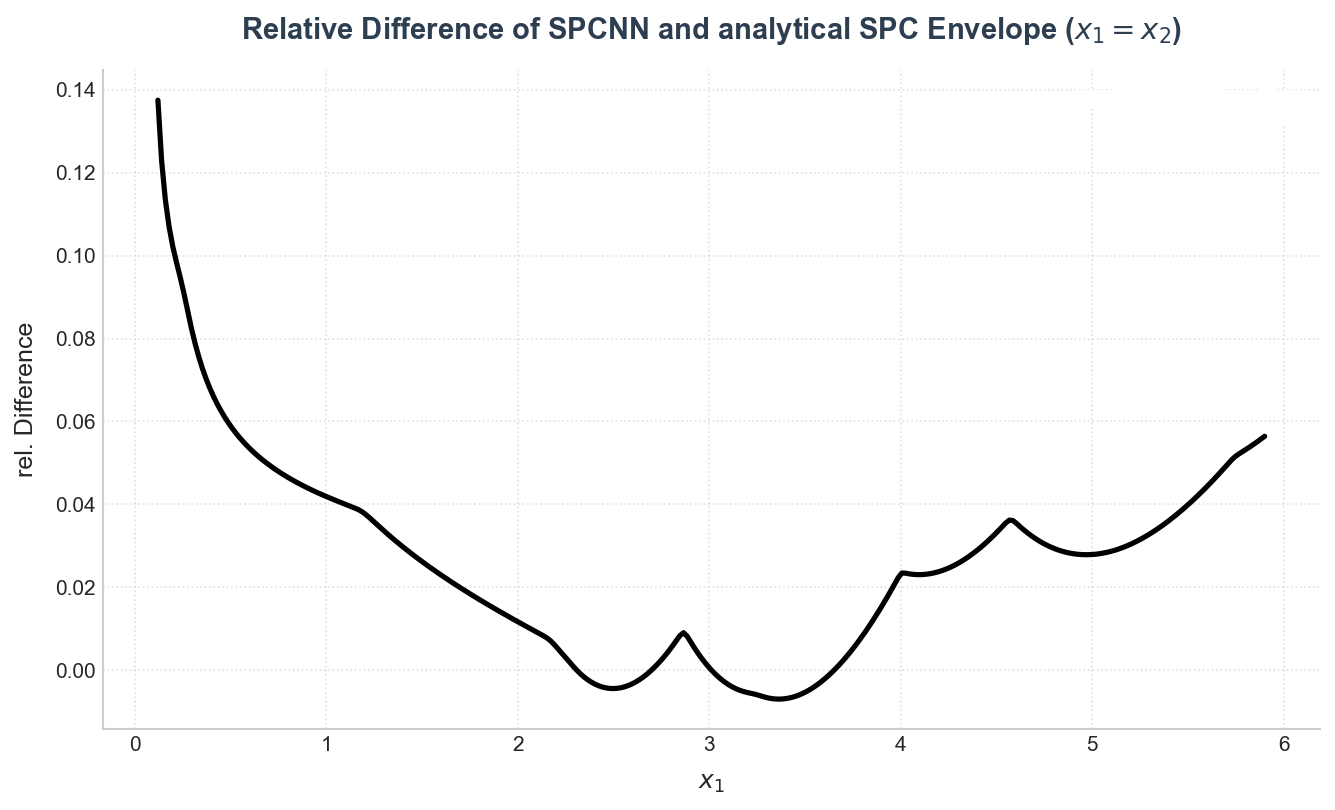}
    \end{minipage}
    
    \caption{Left: total view of the envelope line; Right: relative error}
    \label{Fig: Envelope Error}
\end{figure}

\textbf{Acknowledgements.} 

This work has been funded by the Deutsche Forschungsgemeinschaft (DFG). We gratefully acknowledge the support within the project SPP 2256
(``Variational Methods for Predicting Complex Phenomena in Engineering Structures
and Materials'').

\printbibliography[heading=bibintoc, title={References}]

\end{document}